\documentclass{notices}

\usepackage{aliascnt, hyperref}
\hypersetup{hidelinks}
\usepackage{bbm}
\usepackage{pbox}
\usepackage{arydshln}
\usepackage{pinlabel}
\usepackage{cals}
\usepackage{verbatim}
\usepackage{amsthm}
\usepackage{graphicx, subcaption, wrapfig, framed, multicol, adjustbox}
\usepackage{mathtools}

\newcommand{\gB}{\bold{B}}

\newcommand{\gG}{\bold{G}}

\newcommand{\gP}{\bold{P}}

\newcommand{\bC}{\mathbb{C}}\newcommand{\bD}{\mathbb{D}}
\newcommand{\bF}{\mathbb{F}}

\newcommand{\bQ}{\mathbb{Q}}\newcommand{\bR}{\mathbb{R}}

\newcommand{\bZ}{\mathbb{Z}}
\newcommand{\Z}{\mathbb{Z}}
\newcommand{\Q}{\mathbb{Q}}
\newcommand{\C}{\mathbb{C}}
\newcommand{\R}{\mathbb{R}}

\newcommand{\B}{\text{B}}
\newcommand{\E}{\text{E}}

\newcommand{\FI}{\mathsf{FI}}
\newcommand{\FS}{\mathsf{FS}}
\newcommand{\Mod}{\text{Mod}}
\newcommand{\PMod}{\text{PMod}}

\newcommand{\Diff}{\text{Diff}}
\newcommand{\Aut}{\text{Aut}}
\newcommand{\GL}{\text{GL}}

\usepackage[dvipsnames,svgnames,x11names,hyperref, table]{xcolor}
\usepackage{mathtools,url,verbatim,amssymb,enumerate,stmaryrd}
\usepackage{tikz, tikz-cd}
\usetikzlibrary{matrix,arrows,decorations.pathmorphing}
\usetikzlibrary{arrows,decorations.pathmorphing,backgrounds,fit,positioning,shapes.symbols,chains,calc}

\numberwithin{theoremcounter}{section}
\newaliascnt{theoremauto}{theoremcounter}

\newaliascnt{Defauto}{theoremcounter}

\newaliascnt{exampleauto}{theoremcounter}

\newaliascnt{lemmaauto}{theoremcounter}

\newaliascnt{propositionauto}{theoremcounter}

\newaliascnt{corollaryauto}{theoremcounter}

\newaliascnt{remarkauto}{theoremcounter}

\newaliascnt{notationauto}{theoremcounter}

\newaliascnt{claimauto}{theoremcounter}

\newaliascnt{warningauto}{theoremcounter}

\newaliascnt{questionauto}{theoremcounter}

\newaliascnt{discussionauto}{theoremcounter}

\newaliascnt{computationauto}{theoremcounter}

\newaliascnt{conjectureauto}{theoremcounter}

\newaliascnt{convauto}{theoremcounter}

\newtheorem{theorem}[theoremauto]{Theorem}

\theoremstyle{definition}
\newtheorem{definition}[Defauto]{Definition}

\theoremstyle{remark}

\newtheorem{question}[questionauto]{Question}

\usepackage[vcentermath]{youngtab}
\newcommand{\Y}[1]{\text{\tiny\yng(#1)}}

\title{Stability properties of moduli spaces
}

\author{
  Rita Jim\'enez Rolland
  \affil{
  Rita Jim\'enez Rolland is an Associate Professor at the Institute of Mathematics of UNAM in Oaxaca, Mexico. Her email address is rita@im.unam.mx.
    }
  \and
Jennifer C. H. Wilson
  \affil{Jenny Wilson is an Assistant Professor of mathematics at the University of Michigan in the United States. Her email address is jchw@umich.edu. 
   }
}

\begin{document}

\maketitle


\section{Moduli spaces and stability} 



{\it Moduli spaces} are spaces that parametrize topological or geometric data.    They often appear in families; for example, the configuration spaces of $n$ points in a fixed manifold, 
 the Grassmannians of  linear subspaces of dimension $d$ in $\R^{\infty}$,
and the moduli spaces  $\mathcal{M}_g$ of Riemann surfaces of genus $g$. These families are usually indexed by some geometrically defined quantity, such as the number $n$ of points in a configuration, 
the dimension $d$ of the linear subspaces, or  the genus $g$ of a Riemann surface.  Understanding the topology of these spaces  has been a subject of intense interest for the last 60 years.  

For a  family of moduli spaces  $\{X_n\}_n$ we ask:
\begin{question} \label{QTopology} 
How does the topology of the moduli spaces $X_n$ change as the parameter $n$ changes? 
\end{question}


For many natural examples of moduli spaces $X_n$,  some aspects of the topology get more complicated as the parameter $n$ gets larger. For instance, the dimension of $X_n$  frequently increases with $n$ as well as the number of generators and relations needed to give a presentation of their fundamental groups.  But, maybe surprisingly, there are sometimes features of the moduli spaces that `stabilize' as $n$ increases. In this survey we will describe some forms of stability and some examples of where they arise. 




\subsection{Homology and cohomology} 
Algebraic topology is a branch of mathematics that uses tools from abstract algebra to classify and study topological spaces. By constructing  {\it algebraic invariants}  of topological spaces, we can translate topological problems into (typically easier) algebraic ones. An algebraic invariant of a space is a quantity or algebraic object, such as a group, that is preserved under homeomorphism or homotopy equivalence.  One example is the fundamental group $\pi_1(X,x_0)$ of homotopy classes of loops in a topological space $X$ based at the point $x_0$. Homology and cohomology groups are other examples and are the focus of this article. Their definitions are more subtle than those of homotopy groups like $\pi_1(X,x_0)$, but they are often more computationally tractable and are widely studied.

 Given a topological space $X$ and $k \in \Z_{\geq 0}$, we can associate groups $H_k (X;R)$ and  $H^k(X;R)$,  the $k$th \emph{homology}  and \emph{cohomology groups} (with coefficients in $R$),  where $R$ is a commutative ring such as $\mathbb{Z}$ or  $\mathbb{Q}$.
These algebraic invariants define functors from the category of topological spaces to the category of $R$-modules:    for any
continuous map of topological spaces $f \colon X\rightarrow Y$ there are induced $R$-linear maps $$f_*\colon H_k(X;R)\rightarrow H_k(Y;R)\text{\ \ \ (covariant),}$$
 $$f^*\colon H^k(Y;R)\rightarrow H^k(X;R)\text{\ \ \ (contravariant)}.$$ 
The cohomology groups $H^*(X;R) = \bigoplus_k H^k(X; R)$ in fact have the structure of a graded $R$-algebra with respect to the \emph{cup product} operation. 

The group $H_0(X;\mathbb{Z})$ is  the free abelian group on the path components  of the topological space $X$ and $H^0(X; \mathbb{Z})$ is  its dual. If $X$ is path-connected, $H_1(X;\mathbb{Z})$ is naturally isomorphic to the abelianization of $\pi_1(X, x_0)$ with respect to any basepoint $x_0$, and its elements are certain equivalences classes of (unbased) loops in $X$.


For a topological group $G$ there exists an associated {\it classifying space} $\B G$ for  {\it principal $G$-bundles}. It is constructed as the quotient of a (weakly) contractible space $\E G$ by a proper free action of $G$.  The space $\B G$ is unique up to (weak) homotopy equivalence. 
If $G$ is a discrete group, then  $\B G$ is precisely an {\it Eilenberg-MacLane space $K(G,1)$}, i.e., a  path-connected topological space with $\pi_1(\B G) \cong G$ and trivial higher homotopy groups. For example, up to homotopy equivalence, $\B\Z$ is the circle, $\B \bZ_2$ is the infinite-dimensional real projective space $\mathbb{R}\mathrm{P}^{\infty}$,  
and the Grassmanian of $d$-dimensional linear subspaces in $\mathbb{R}^{\infty}$ is $\B \GL_d(\mathbb{R})$. 

Some motivation to study the cohomology of $\B G$: its cohomology classes define \emph{characteristic classes}  of principal $G$-bundles, invariants that measure the `twistedness' of the bundle. For instance the cohomology algebra $H^*( \B\GL_d(\mathbb{R});\mathbb{Z})$ can be described in terms of Pontryagin and  Stiefel--Whitney classes. 


With $\B G$ we can define  the {\it group homology}  and {\it group cohomology} of  a discrete group $G$ by $$H_k(G;R) \vcentcolon= H_k(\B G;R),\ \ H^k(G; R) \vcentcolon= H^k(\B G; R).$$

We can refine our   \autoref{QTopology} to the following:

\begin{question} \label{QHomology} Given  family $\{X_n\}_{n}$ of moduli spaces or discrete groups,  how do the homology and cohomology groups of  the $n$th space in the sequence change as the parameter $n$ increases?
\end{question}

In this  article we discuss \autoref{QHomology} with a particular focus on the families of configuration spaces and braid groups.  For further reading 
we recommend R. Cohen's survey \cite{CohenSurvey}  on stability of moduli spaces.

\subsection{Homological stability}


\begin{definition} 
A sequence of spaces  or groups with maps $$X_0 \xrightarrow{s_{0}} \ldots \xrightarrow{s_{n-2}} X_{n-1}\xrightarrow{s_{n-1}} X_n\xrightarrow{s_n} X_{n+1}\xrightarrow{s_{n+1}}\ldots$$ satisfies {\it homological stability} if, for each $k$, the  induced map in degree-$k$ homology $$(s_n)_*\colon H_k(X_n;\bZ)\rightarrow H_k(X_{n+1};\bZ)$$ is an isomorphism for all $n \geq N_k$ for some  stability threshold $N_k \in \Z$ depending on $k$. The maps $s_n$ are sometimes called \emph{stabilization maps} and the set $\{ (n,k) \in \Z^2 \; | \; n \geq N_k\}$  is the {\it stable range}. 
\end{definition} 

 If the maps $s_n\colon  X_n \to X_{n+1}$ are inclusions we define  $X_{\infty} := \bigcup_{n\geq 1} X_n$ to be the \emph{stable} group  or space. Under mild assumptions, if $\{X_n\}_n$ satisfies homological stability, then $$H_k(X_{\infty};\bZ)\cong H_k(X_n;\bZ) \text{ \ \  for \ \  }n \geq N_k.$$  We call the groups $H_k(X_{\infty};\bZ)$ the {\it stable homology}.   

 





\section{An example: configuration spaces and the braid groups}

\subsection{A primer on configuration spaces} 
\begin{definition} 
Let $M$ be a topological space, such as a graph or a manifold.  The \emph{(ordered) configuration space $F_n(M)$ of $n$ particles on $M$} is the space 
$$ F_n(M) = \{ (x_1,  \ldots, x_n) \in M^n \; | \; x_1, \ldots, x_n \text{ distinct} \},$$ 
topologized as a subspace of $M^n$. Notably, $F_0(M)$ is a point and $F_1(M)=M$.  
\end{definition} Configuration spaces have a long history of study in connection to topics as  broad-ranging as homotopy groups of spheres and robotic motion planning.

\begin{figure}[h!]
\raggedleft
\centering
\raisebox{.8em}{$F_2([0,1]) =$}   \includegraphics[scale=.8]{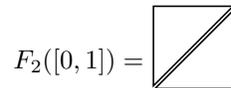} 
\caption{The space $F_2([0,1])$ is obtained by deleting the diagonal from the square $[0,1]^2$.} 
\label{F2I} 
\end{figure}

One way to conceptualize the configuration space $F_n(M)$ is as the complement of the union of subspaces of $M^n$ defined by equations of the form $x_i = x_j$. In other words, we can construct $F_n(M)$ by deleting the ``fat diagonal" of $M^n$, consisting of all $n$-tuples in $M^n$ where two or more components coincide. In the simplest case, when $n=2$ and $M$ is the interval $[0,1]$,  we see that $F_2([0,1])$ consists of two contractible components, as in \autoref{F2I}. 

\begin{figure}[h!]
\raggedleft
\centering
 \includegraphics[scale=.25]{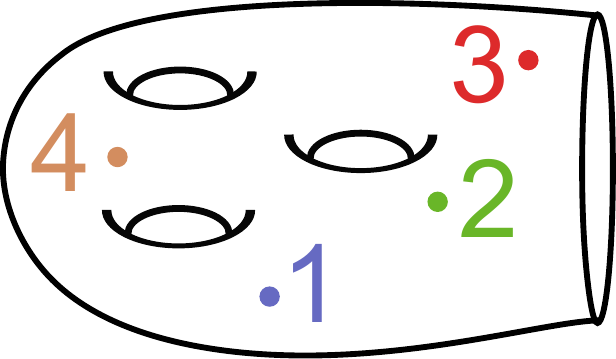}  \raisebox{1em}{  $\in F_4(\Sigma)$}  
 \caption{ A point in the ordered configuration space of an open surface $\Sigma$.} 
\label{PointinF4M} 
\end{figure}

Another way we can conceptualize $F_n(M)$ is as the space of embeddings of the discrete set $\{1, 2, \ldots, n\}$ into $M$, appropriately topologized. We may visualize a point in $F_n(M)$ by labelling $n$ points in $M$, as in \autoref{PointinF4M}. 

\begin{figure}[h!]
\raggedleft
\centering
\raisebox{.8em}{$F_2([0,1]) =$}    \includegraphics[scale=.8]{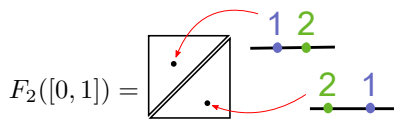} 
\caption{ The path components of $F_2([0,1])$.} 
\label{F2IComponents} 
\end{figure}

From this perspective, we may reinterpret the path components of $F_2([0,1])$: one component consists of all configurations where particle 1 is to the left of particle 2, and one component has particle 1 on the right. See  \autoref{F2IComponents}. 

\begin{figure}[h!]
\raggedleft
\centering
\includegraphics[scale=.25]{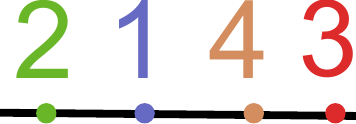}  \; \;  \raisebox{.7em}{  $\in F_4([0,1]) $} 
\caption{A point in $F_4([0,1])$ in the path component indexed by the permutation $2143$ in $S_4$. } 
\label{F4I} 
\end{figure}

Any path through $[0,1]^2$ that interchanges the relative positions of the two particles must involve a `collision' of particles, and hence exit the configuration space $F_2([0,1]) \subseteq [0,1]^2$. We encourage the reader to verify that, in general, the configuration space $F_n([0,1])$ is the union of $n!$ contractible path components, indexed by elements of the symmetric group $S_n$. See \autoref{F4I}. 

In contrast, if $M$ is a connected manifold of dimension $2$ or more, then $F_n(M)$ is path-connected: given any two configurations, we can construct a path through $M^n$ from one configuration to the other without  any `collisions' of particles. In this case $H_0(F_n(M);\bZ) \cong \Z$ for all $n \geq 0$, and this is our first glimpse of stability in these spaces as $n \to \infty$. 

\begin{figure}[h!]
\raggedleft
\centering
\begin{tikzpicture}
  \matrix (m) [matrix of math nodes,row sep=4.5em,column sep=4em,minimum width=2em]
  {  F_n(M) \\ C_n(M) := F_n(M) / S_n \\};
  \path[-stealth]
    (m-1-1) edge (m-2-1);
\end{tikzpicture}
\includegraphics[scale=.19]{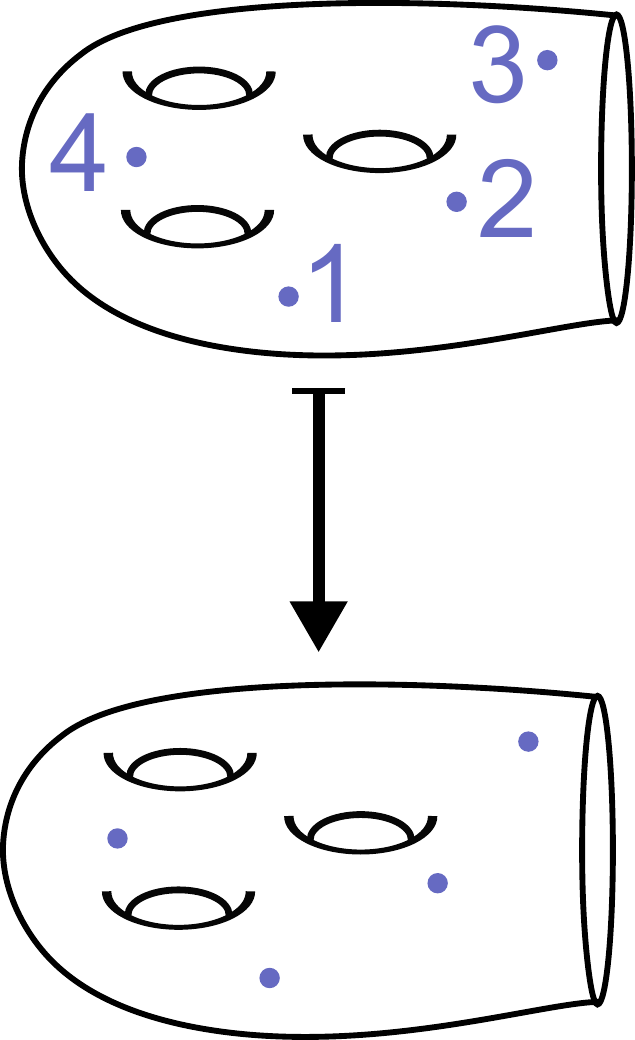} 
\caption{The quotient map $F_n(M) \to C_n(M)$ } 
\label{CnMQuotient} 
\end{figure}

For any space $M$, the symmetric group $S_n$ acts freely on $F_n(M)$ by permuting the coordinates of an $n$-tuple $(x_1, \ldots, x_n)$, equivalently, by permuting the labels on a configuration as in  \autoref{PointinF4M}.  The orbit space $C_n(M) = F_n(M) /S_n$ is the \emph{(unordered) configuration space of $n$ particles on $M$}. This is the space of all $n$-element subsets of $M$, topologized as the quotient of $F_n(M)$. The reader may verify that the quotient map (illustrated in  \autoref{CnMQuotient}) is a regular $S_n$-covering space map. In particular, by covering space theory, the quotient map $F_n(M) \to C_n(M)$ induces an injective map 
on fundamental groups.

In the case that $M$ is the complex plane $\C$, we can identify  $C_n(\C)$ with the space of monic degree-$n$ polynomials over $\C$ with distinct roots, by mapping a configuration $\{z_1, \ldots, z_n\}$ to the polynomial $p(x)=(x-z_1)\cdots(x-z_n)$. For this reason the topology of $C_n(\C)$ has deep connections to classical problems about finding roots of polynomials. 

 We will address \autoref{QHomology} for the families $\{ C_n(M)\}_n$ and $\{ F_n(M)\}_n$, but  we first specialize to the case when $M=\bC$.  Although the spaces $C_n(\bC)$ and $F_n(\bC)$ are path-connected, in contrast to the configuration spaces of $M=[0,1]$, they have rich topological structures: they are classifying spaces for the braid groups and the pure braid groups, respectively, 
 which we now introduce.  

 \subsection{A primer on the braid groups}   \label{SectionBraid} 
\begin{figure}[h!]
\raggedleft
\centering
\includegraphics[scale=.6]{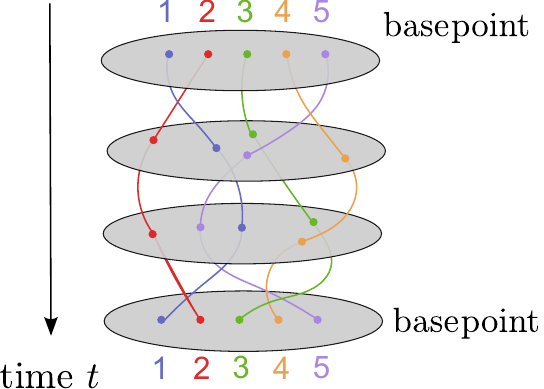} 
\caption{A visualization of a loop $\gamma(t)$ in $F_5(\bC)$  representing an element of $\pi_1(F_5(\bC)) \cong \gP_5$. } 
\label{BraidAsLoop} 
\end{figure}

Since $F_n(\C)$ is path-connected, as an abstract group its fundamental group is independent of choice of basepoint. For path-connected spaces, we sometimes drop the basepoint from the notation for $\pi_1$. 
 \begin{definition}The fundamental group $\pi_1(C_n(\bC))$  is called the  {\it braid group} $\gB_n$ and $\pi_1(F_n(\bC))$ is the {\it pure braid group } $\gP_n$.  
\end{definition}
We can understand $\pi_1(F_n(\bC))$ as follows. Choose a basepoint configuration $(z_1, \ldots, z_n)$ in  $F_n(\bC)$, and then we may visualize a loop as a `movie' where the $n$ particles continuously move around $\bC$, eventually returning pointwise to their starting positions. 
If we represent time by a third spacial dimension, as shown in  \autoref{BraidAsLoop}, we can view the particles as tracing out a braid. Note that, up to homeomorphism, we may view $F_n(\C)$ as the configuration space of the open 2-disk.

Loops in  $C_n(\bC)$ are similar, with the crucial distinction that the $n$ particles are unlabelled and indistinguishable, and so need only return set-wise to their basepoint configuration.

\begin{figure}[h!]
\centering
\includegraphics[scale=.5]{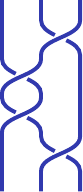} 
\caption{A braid on 3 strands} 
\label{Braid-3Strands} 
\end{figure}

\begin{figure}[h!]
\raggedleft
\centering
\includegraphics[scale=.5]{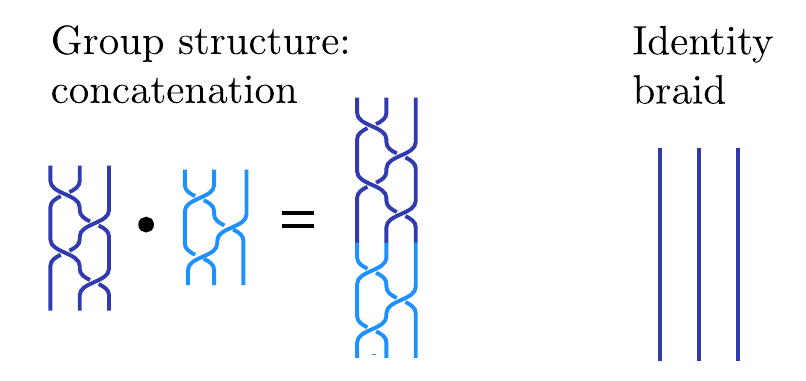} 
\caption{The group structure on $\gB_n$} 
\label{Braid-GroupOperation&Identity} 
\end{figure}

It is traditional to represent elements of the group $\gB_n$ and its subgroup $\gP_n$ by equivalence classes of \emph{braid diagrams}, as illustrated in  \autoref{Braid-3Strands}. 
These braid diagrams depict $n$ strings (called \emph{strands}) in Euclidean 3-space, anchored at their tops at $n$ distinguished points in a horizontal plane, and anchored at their bottoms at the same $n$ points in a parallel plane. The strands may move in space but may not double back or pass through each other. The group operation is concatenation, as in  \autoref{Braid-GroupOperation&Identity}. 
%


The braid groups were defined rigorously by Artin in 1925 \cite{Artin}, but the roots of this notion appeared already in the work of Hurwitz,  Firckle, and Klein in  1890's, and of Vandermonde in 1771. This topological interpretation  of braid groups as the fundamental groups of configuration spaces was formalized in 1962 by Fox and Neuwirth \cite{FoxN}.

\begin{figure}[h!]
\raggedleft
\centering
\includegraphics[scale=.6]{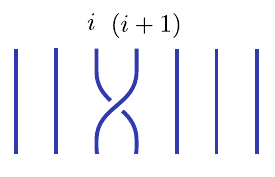} 
\caption{Artin's generator $\sigma_i$ for $\gB_n$} 
\label{StandardGenerator} 
\end{figure}

Artin established presentations for the braid group and the pure braid group. His presentation for $\gB_n$,
\begin{equation*}
 {\small \gB_n \cong  \left\langle \sigma_1,\sigma_2\dots, \sigma_{n-1}  \biggm|
\begin{array}{c}   \sigma_i\sigma_j = \sigma_j\sigma_i   \ \text{if}  \  |i-j| \geq 2 \\
\sigma_i\sigma_{i+1}\sigma_i  = \sigma_{i+1}\sigma_i\sigma_{i+1} 
\end{array}
\right\rangle },
\end{equation*}
uses $(n-1)$ generators  $\sigma_i$ corresponding to half-twists of adjacent strands, as in  \autoref{StandardGenerator}. 
%
%
%

\begin{figure}[h!]
\raggedleft
\centering
\includegraphics[scale=.8]{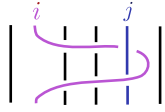} 
\caption{Artin's generator $T_{ij}=T_{ji}$ for $\gP_n$} 
\label{PureBraidGenerator} 
\end{figure}

Artin also gave a finite presentation for $\gP_n$. We will not state it in full, but comment that there are ${n \choose 2}$ generators $T_{ij}$, ($i \neq j$, $i,j \in \{1, 2, \ldots, n\})$ corresponding to full twists of each pair of strands, as in  \autoref{PureBraidGenerator}. 

\begin{figure}[h!]
\raggedleft
\centering
\includegraphics[scale=.7]{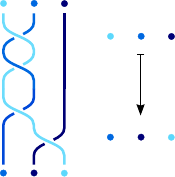} 
\caption{The quotient map $\gB_n \to S_n$} 
\label{BntoSn} 
\end{figure}

Corresponding to the regular covering space map $F_n(\C) \to C_n(\C)$ of  \autoref{CnMQuotient}, there is a short exact sequence of groups 
$$ 1 \to \gP_n \to \gB_n \to S_n \to 1.$$
The quotient map $ \gB_n \to S_n$, shown in  \autoref{BntoSn}, takes a braid, forgets the $n$ strands, and simply records the permutation induced on their endpoints. The generator $\sigma_i$ maps to the simple transposition $(i \; i+1)$. The kernel is those braids that induce the trivial permutation, i.e.,  the pure braid group. 

\subsection{Homological stability for the braid groups} 

Arnold \cite{ArnoldStability} calculated some homology groups of $\gB_n$ in low degree (\autoref{LowHomology}). 

\begin{table}[h]
\begin{center} {\footnotesize
\begin{tabular}{|ccc|c|c|c|c|c|c|} \hline
&&&  &  &  &  &  &  \\[-.5em]
&& $k$ & 0 & 1 & 2 & 3 & 4 & 5 \\[-.5em]
&$n$ &&  &  &  &  &  &  \\ \hline
&0 && $ \cellcolor{blue!15} \bZ$ &  &  &  &  &   \\
&1 && $ \cellcolor{blue!15} \bZ$ &  &  &  &  &   \\
&2 &&  \cellcolor{blue!15}$\bZ$ & $ \cellcolor{blue!15}{\bZ}$ &  &  &  &   \\
&3 && \cellcolor{blue!15}$\bZ$ & $ \cellcolor{blue!15}\bZ$ &  &  &  &   \\
&4 &&  \cellcolor{blue!15}$\bZ$ & $ \cellcolor{blue!15}\bZ$ &$ \cellcolor{blue!15}\bZ_2$ &  &  &   \\
&5 &&  \cellcolor{blue!15}$\bZ$ & $ \cellcolor{blue!15}\bZ$ &$ \cellcolor{blue!15}\bZ_2$  &  &  &   \\
&6 &&  \cellcolor{blue!15}$\bZ$ & $ \cellcolor{blue!15}\bZ$ &$ \cellcolor{blue!15} \bZ_2$ &$ \cellcolor{blue!15}{ \bZ_2}$  &$\bZ_3$  &    \\
&7 &&  \cellcolor{blue!15}$\bZ$ & $ \cellcolor{blue!15}\bZ$ &$ \cellcolor{blue!15}\bZ_2$  &$ \cellcolor{blue!15}\bZ_2$ &$\bZ_3$  &     \\
&8 &&  \cellcolor{blue!15}$\bZ$ & $ \cellcolor{blue!15}\bZ$ &$ \cellcolor{blue!15}\bZ_2$  &$ \cellcolor{blue!15}\bZ_2$ &$ \cellcolor{blue!15}{\bZ_6}$  &$\bZ_3$    \\
&9 && \cellcolor{blue!15}$\bZ$ & \cellcolor{blue!15}$\bZ$ &\cellcolor{blue!15}$\bZ_2$  &\cellcolor{blue!15}$\bZ_2$ &\cellcolor{blue!15}$\bZ_6$  &$\bZ_3$     \\ 
\hline
\end{tabular} }
\end{center}
\caption{\label{LowHomology} The homology groups $H_k(\gB_n;\bZ)$. Empty spaces are zero groups. Stable groups are shaded.}
\end{table}
The $k=0$ column follows from the fact that 
$C_n(\bR^2)$ is path-connected  and the $k=1$ column can be obtained by abelianizing  Artin's presentation of $\gB_n$.  
 Even the low-degree calculations in  \autoref{LowHomology} suggest a pattern: the  homology of $\gB_n$ in a fixed degree $k$ becomes independent of $n$ as $n$ increases.

\begin{figure}[h!]
\raggedleft
\centering
\includegraphics[scale=.5]{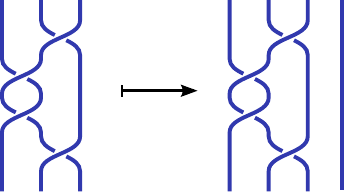} 
\caption{The stabilization map $s_3\colon \gB_3\hookrightarrow \gB_{4}$ } 
\label{Braid-stabilization} 
\end{figure}

Arnold \cite{ArnoldStability} proved the following stability result, in terms of the stabilization map $s_n\colon \gB_n\hookrightarrow \gB_{n+1}$  defined by 
 adding an unbraided $(n+1)^{st}$ strand as in \autoref{Braid-stabilization}.  

\begin{theorem}[Arnold \cite{ArnoldStability}] \label{HomStabilityBraids} 
For each $k\geq 0$, the induced map   $$(s_n)_*\colon H_k(\gB_n;\bZ)\rightarrow H_k(\gB_{n+1};\bZ)$$ is an isomorphism for 
 $n\geq 2k$.
\end{theorem}

The family $\{C_n(\bC)\}_n$ therefore satisfies homological stability. 
Arnold \cite{ArnoldStability} in fact proved the result for cohomology, and \autoref{HomStabilityBraids} follows from the universal coefficients theorem.

May \cite{May} and Segal \cite{Segal} proved that the  stable braid group $\gB_{\infty}$ has the same homology as the path component of the trivial loop in the double loop space $\Omega^2 S^2$. 
F. Cohen \cite{CohenBraids}  
and Va\u{\i}n\u{s}te\u{\i}n \cite{VBraids}.
computed the cohomology ring with coefficients in $\mathbb{F}_p$ ($p$ an odd prime), and described $H^k(\gB_n; \Z)$ in terms of the groups $H^{k-1}(\gB_n; \mathbb{F}_p)$ ($p \text{ prime})$ for $k \geq 2$. 
 \subsection{Homological stability  for configuration spaces}\label{HomStabConf}

\begin{figure}[h!]
\raggedleft
\centering
\includegraphics[scale=.6]{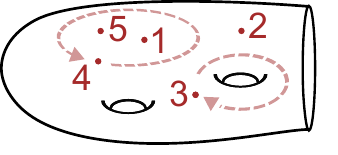} 
\caption{A class in $H_2(F_5(M))$} 
\label{HomologyDemo} 
\end{figure}

For a $d$-manifold $M$, it is possible to visualize homology classes in $F_n(M)$ and $C_n(M)$ concretely. Consider \autoref{HomologyDemo}.
This figure shows a $2$-parameter family of configurations in $F_n(M)$; in fact (because the two loops do not intersect) it shows an embedded torus $S^1 \times S^1 \hookrightarrow F_5(M)$. Thus, up to sign, this figure represents an element of $H_2(F_5(M))$. In a sense, the loop traced out by particle $3$ arises from the homology of the surface $M$, and the loop traced out by particle $4$ arises from the homology of $F_n(\R^d)$. From the homology of $M$ and $F_n(\R^d)$, it is possible to generate lots of examples of homology classes in $F_n(M)$. The problem of understanding additive relations among these classes, however, is subtle, and the groups $H_k(F_n(M);\bZ)$ are unknown in most cases. 

When $M$ is (punctured) Euclidean space, the (co)homology groups of $F_n(M)$ were computed by  Arnold and Cohen. However, even in the case that $M$ is a genus-$g$ surface, we currently do not know the Betti numbers $\beta_k = \mathrm{rank}(H_k(F_n(M);\bZ))$. Recently Pagaria \cite[Corollary 2.9]{Pagaria-Ordered} computed the asymptotic growth rate of the Betti numbers in the case $M$ is a torus. In the case of {\it unordered} configuration spaces, in 2016 Drummond-Cole and Knudsen \cite{DCK} computed the Betti numbers of $C_n(M)$ for $M$ a surface of finite type.

Even though the (co)homology groups of configurations spaces remain largely mysterious, the tools of homological stability give us a different approach to understanding their structure. 

\begin{figure}[h!]
\raggedleft
\centering
\includegraphics[scale=.6]{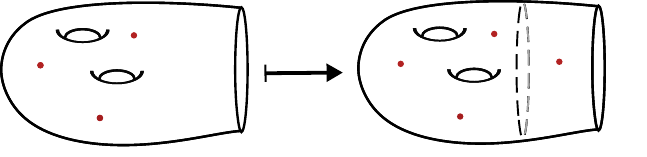} 
\caption{Stabilization map  $s_3\colon C_3(M)\rightarrow C_4(M)$} 
\label{Surfaces} 
\end{figure}

\autoref{HomStabilityBraids} on stability for braid groups raises the question of whether the unordered configurations spaces $\{C_n(M)\}_n$ satisfy homological stability for a larger class of topological spaces $M$.  Let $M$ be a connected manifold. To generalize \autoref{HomStabilityBraids} we must define stabilization maps 
\begin{align*}
C_n(M) & \longrightarrow C_{n+1}(M) \\ 
\{x_1, \ldots, x_n\} & \longmapsto \{x_1, \ldots, x_n, x_{n+1}\}
\end{align*} 
Unfortunately, in general there is no way to choose a distinct particle $x_{n+1}$ continuously in the inputs $\{x_1, \ldots, x_n\}$, and no continuous map of this form exists. To define the stabilization maps, we must assume extra structure on $M$, for example, assume that $M$ is the interior of a manifold with nonempty boundary. Then, if we choose a boundary component, it is possible to define the stabilization map $s_n\colon  C_n(M) \to C_{n+1}(M)$ by placing the new particle in a sufficiently small collar neighbourhood of the boundary component. This procedure (illustrated in \autoref{Surfaces}) is informally described as `adding a particle at infinity.' 

In the 1970's McDuff 
proved that the sequence $\{C_n(M)\}_n$ satisfies homological stability and Segal 
gave explicit stable ranges.  
\begin{theorem}[McDuff \cite{McDuff};  Segal \cite{SegalRational}]  Let $M$ be the interior of a compact connected manifold with nonempty boundary. For each $k \geq 0$ the maps
$$ (s_n)_*\colon   H_k(C_n(M);\bZ) \longrightarrow  H_k(C_{n+1}(M);\bZ) $$ are isomorphisms for $n \geq 2k$. 
\end{theorem} 
    
\begin{figure}[h!]
\raggedleft
\centering
\includegraphics[scale=.75]{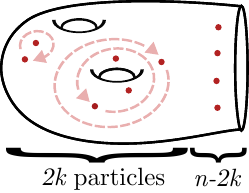} 
\caption{A homology class after stabilizing by the addition of $n-2k$ particles. } 
\label{ImageOfStabilization} 
\end{figure}

 Concretely, this theorem states that degree-$k$ homology classes arise from subconfigurations on at most $2k$ particles. Heuristically, these homology classes have the form of \autoref{ImageOfStabilization}. 

Moreover,  McDuff related the homology of the stable space $C_{\infty}(M)$  to the homology of  $\Gamma(M)$,  the space of compactly-supported smooth sections of the bundle over $M$ obtained by taking the fibrewise one-point compactification of the tangent bundle of $M$. 





  


\section{Other stable families} 

We briefly describe some other significant families satisfying (co)homological stability.  \\[-.5em]

\noindent{\bf Symmetric groups.}  Nakaoka  \cite{Nakaoka}  proved  that the symmetric groups $\{S_n\}_n$ satisfy homological stability with  respect to the inclusions 
$S_n\hookrightarrow S_{n+1}$ (see also \cite{Kerz}). The Barratt--Priddy--Quillen theorem \cite{BarrattPriddy} states that the infinite symmetric group $S_{\infty} = \bigcup_n S_n$ has the
same homology of $\Omega_0^{\infty} S^{\infty}$, the path-connected component of the identity in the infinite loop space $\Omega^{\infty} S^{\infty}$. \\[-.5em] 


\noindent{\bf General linear groups.}  Let $R$ be a ring. Consider  the sequence of general linear groups $\{\GL_n(R)\}_n$
with the inclusions
$\GL_n(R)\hookrightarrow \GL_{n+1}(R)$ given by  $$ A\mapsto\begin{bmatrix}A & 0 \\0 & 1 \end{bmatrix}.$$
In the 1970's Quillen  studied the homology of these groups when $R$ is a finite field $\bF_q$ of characteristic $p$ in his seminal work \cite{Quillen} on the $K$-theory of finite fields. He computes $H^*(\GL_n(\bF_q); \bF_{\ell})$ for prime $\ell \neq p$  and determines a vanishing range for $\ell =p$. 



Charney \cite{Charney} proved homological stability when  $R$  is a Dedekind domain. Van der Kallen \cite{vanderKallen}, building on work of Maazen \cite{Maazen}, proved the case that $R$ is an associative ring satisfying Bass's ``stable rank condition"; this arguably includes any naturally-arising ring. 

These results are part of a large stability literature on classical groups that warrants its own survey. It includes work by Betley, Cathelineau, Charney, Collinet,  Dwyer, Essert, Friedlander, Friedrich  Galatius,  Guin,  Hutchinson, Kupers,  Miller,  Mirzaii,  Nesterenko,  Panin, Patzt,  Randal-Williams, Sah,   Sprehn, Suslin,   Tao,   Vaser\v{s}te\u{\i}n,   Vogtmann, and Wahl \cite{Vaser, Friedlander, Charney,  Dwyer, Vogt-GL, Suslin,  Sah, Guin, Panin,  Nester,  Betley, Mirzaii, Mirzaii2, Catheline,  Tao, Collinet, Essert,  Fried, GKRW-GL,  KMP-GL,  SW  }, among others. Homological stability is known to hold for special linear groups, orthogonal groups, unitary groups, and other families of classical groups. There is ongoing work to study (co)homology with twisted coefficients, and sharpen the stable ranges.   \\[-.5em]





\begin{figure}[h!]
\raggedleft
\centering
\includegraphics[width=.95\columnwidth]{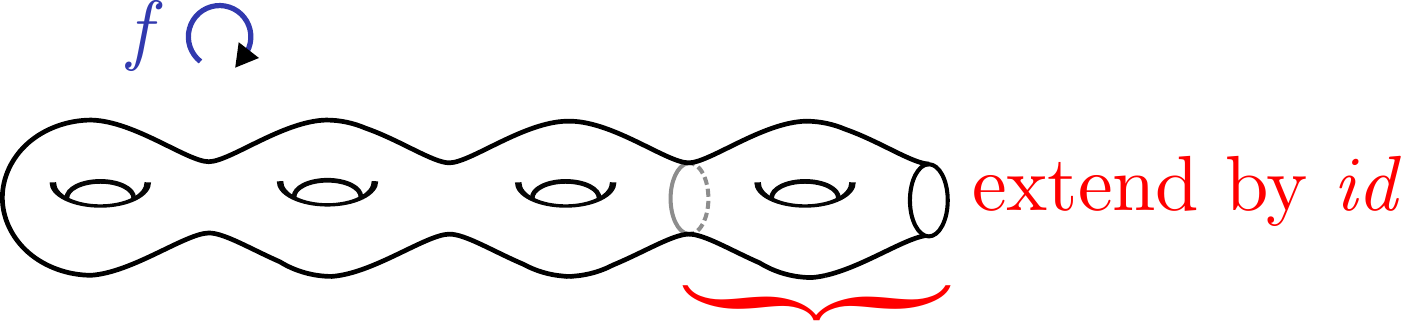} 
\caption{The map  $\Mod(\Sigma_{3,1})\to \Mod(\Sigma_{4,1})$  is induced by the inclusion $\Sigma_{3,1}\hookrightarrow \Sigma_{4,1}$} 
\label{MCGinclusion} 
\end{figure}

\noindent{\bf Mapping class groups and moduli space of Riemann surfaces.}  Let $\Sigma_{g,1}$ be an oriented surface of genus $g$ with one boundary component and let the {\it mapping class group} $$\Mod(\Sigma_{g,1}):=\pi_0(\Diff^+(\Sigma_{g,1}\text{ rel }\partial))$$ be the group of isotopy classes of diffeomorphisms of $\Sigma_{g,1}$ 
fixing a collar neighbourhood of the boundary. There is a map $t_g\colon \Mod(\Sigma_{g,1})\hookrightarrow \Mod(\Sigma_{g+1,1})$  induced by the inclusion $\Sigma_{g,1}\hookrightarrow \Sigma_{g+1,1}$ by extending a diffeomorphism by the identity on the complement $\Sigma_{g+1,1}\setminus \Sigma_{g,1}$, as in \autoref{MCGinclusion}.

There is also a map $cap\colon \Mod(\Sigma_{g,1})\rightarrow \Mod(\Sigma_{g})$ induced by gluing a disk on the boundary component of $\Sigma_{g,1}$.  Harer \cite{Harer} proved that the sequence $\{ \Mod_{g,1}\}_{g}$ satisfies homological stability with respect to the inclusions $t_g$ and that for large $g$ the map $cap$ induces isomorphisms on homology. 
The proof and the stable ranges have been improved by  work of  Ivanov, Boldsen and others \cite{Ivanov, Boldsen}.  Madsen and Weiss   \cite{MadsenWeiss}  computed the stable homology by identifying  the homology of mapping class groups, in the stable range, with the homology of a certain infinite loop space.  
 
The rational homology of the  mapping class group $\Mod(\Sigma_g)$ is  the same as that of  the {\it  moduli space $\mathcal{M}_g$ of Riemann surfaces}  of genus $g\geq 2$. This moduli space  parametrizes:
 \begin{itemize} 
 \item isometry classes of hyperbolic structures on $\Sigma_g$,
\item conformal classes of Riemannian metrics on  $\Sigma_g$,
\item biholomorphism classes of complex structures on  the surface $\Sigma_g$,
\item isomorphism classes of smooth algebraic curves homeomorphic to $\Sigma_g$.
  \end{itemize}

One consequence of Harer's stability theorem and the Madsen--Weiss's theorem is their proof of {\it Mumford's conjecture} \cite{Mumford}: the rational cohomology of $\mathcal{M}_g$ is a polynomial algebra on  generators $\kappa_i$ of degree $2i$, the so-called  {\it Mumford--Morita--Miller classes}, in a stable range depending on $g$. See Tillman's survey \cite{Tillman-MumfordSurvey}
and Wahl's survey \cite{WahlSurvey}.
  

 
Homological stability was established for mapping class groups of non-orientable surfaces by Wahl \cite{Wahl}, for mapping class groups of some $3$-manifolds by Hatcher--Wahl \cite{HW} and  framed, Spin, and Pin mapping class groups by Randal-Williams \cite{SpinMCG}.\\[-.5em]


\noindent{\bf  Automorphism groups of free groups.}  Let $F_n$ denote the free group of rank $n$.   Hatcher and Vogtmann  \cite{HV} proved  that the sequence  $\{ \Aut(F_n)\}_n$ 
is homologically stable with respect to inclusions $ \Aut(F_n)\hookrightarrow \Aut(F_{n+1})$.  Galatius \cite{GalatiusFree}  computed the stable homology by proving that $H_*(\Aut(F_{\infty});\bZ)\cong H_*(\Omega_0^{\infty} S^{\infty};\bZ)\cong H_*(S_{\infty};\bZ)$. In particular, for $n>2k+1$,
$$H_k(\Aut(F_n);\bQ)\cong H_k(\Aut(F_{\infty});\bQ)=0.$$

\noindent{\bf Moduli spaces of high-dimensional manifolds.}  Let $M$ be a smooth compact manifold. The moduli space $\mathcal{M}(M)$ of manifolds of type $M$ is the classifying space $\B\Diff (M \text{\ rel\ } \partial)$. In the last few years Galatius and Randal-Williams \cite{HighDim1} proved homological stability for $\mathcal{M}(M)$ for simply connected manifolds $M$ of dimension $2d>4$, with respect to the $n$-fold connected sum with $S^d \times S^d$. This generalizes Harer's result to higher-dimensional manifolds. They also obtained a generalized Madsen--Weiss theorem for simply connected  manifolds of dimension  $2d>4$ \cite{HighDim2}. Homological stability  with respect to connected sum with $S^p\times S^q$, for $p<q<2p-2$ was obtained by Perlmutter \cite{Perlmutter}. 


\section{A proof strategy}

There is a well-established strategy for proving homological stability that traces back to unpublished work by Quillen in the 1970's  \cite{QuillenUnpublished}.  We describe a simplified version of  Quillen's argument  for a family of discrete groups with inclusions. 

Recall that a \emph{$p$-simplex} $\Delta^p$ is a $p$-dimensional polytope defined as the convex hull of $(p+1)$ points in $\R^{p}$ in general position, called its \emph{vertices}. For example, a $0$-simplex is a point, a $1$-simplex is a closed line segment, and a 2-simplex is triangle.  A \emph{face} of a simplex is the convex hull of a subset of its vertices. A map $f\colon  \Delta^p \to \Delta^q$  is \emph{simplicial} if it maps vertices to vertices, and takes the form
$$ f\colon  \sum_{i=0}^p t_i v_i \mapsto \sum_{i=0}^p t_i f(v_i) $$ with $v_0, \ldots, v_p$ the vertices of $\Delta^p$ and $0 \leq t_i \leq 1$, $ \sum_i t_i =1$. 

 A \emph{triangulation} of a topological space $W$ is a decomposition of $W$ as a union of simplices, such that the intersection $\sigma \cap \tau$ of any pair of simplices $\sigma, \tau$ in $W$ is either empty or equal to a single common face of  $\sigma$ and $\tau$.  A triangulated space is called a \emph{simplicial complex}.  A map $f$ of simplicial complexes is \emph{simplicial} if it maps simplices to simplices and its restriction to each simplex is simplicial. 
 

A simplicial complex $W$ is called \emph{$(-1)$-connected} if it is nonempty, \emph{$0$-connected} if it is path-connected, and \emph{$1$-connected} if it is simply connected. More generally, a nonempty simplicial complex $W$  is called \emph{ $d$-connected} if its homotopy groups $\pi_i(W)$ vanish for all $0 \leq i \leq d$. By the Hurewicz theorem, $W$ is $d$-connected ($ d\geq 2$) if and only if $W$ is simply connected and $H_i(X) =0$ for all $2 \leq i \leq d$. 

With this terminology, we can now describe Quillen's argument. The following formulation of \autoref{QuillenStabillityTheorem} is due to Hatcher--Wahl \cite[Theorem 5.1]{HW}.


\begin{theorem}[\bf Quillen's argument for homological stability] \label{QuillenStabillityTheorem} Let  
$0 \hookrightarrow G_1  \hookrightarrow \ldots \hookrightarrow G_{n} \hookrightarrow\ldots$ be a sequence of discrete groups.  For each $n$ let $W_n$ be a simplicial complex with a simplicial action of $G_n$ satisfying the following properties:

\begin{enumerate}[(i)]
\item The simplicial complexes $W_n$ are $\left(\frac{n-2}{2}\right)$-connected. 

\item For each $p\geq 0$, the group $G_n$ acts transitively on the set of $p$-simplices.

\item For each simplex $\sigma_p$ in $W_n$, the stabilizer $stab(\sigma_p)$  fixes $\sigma_p$ pointwise.

\item The stabilizer $stab(\sigma_p)$ of a $p$-simplex $\sigma_p$ is conjugate in $G_n$ to the subgroup $G_{n-p-1}\subseteq G_n$. (By convention $G_n=0$ if $n\leq 0$.)

\item For each edge $[v_0, v_1]$ in $W_n$, there exists  $g \in G_n$ such that $g \cdot v_0 = v_1$ and $g$ commutes with all elements of $G_n$ that fix $[v_0,v_1]$ pointwise. 

\end{enumerate}

Then the sequence $\{G_n\}_n$ is homologically stable. Specifically, the inclusion $G_n \hookrightarrow G_{n+1}$ induces an isomorphism on degree-$k$ homology for $n \geq 2k+1$ and a surjection for $n=2k$. 
\end{theorem}


Theorem \ref{QuillenStabillityTheorem} follows from a formal algebraic argument involving a sequence of spectral sequences associated to the complexes $W_n$.  
We remark, for the readers familiar with spectral sequences, that for each $n$ we obtain a homology spectral sequence by using  $W_n\times_{G_n} \E G_n$ to build an approximation to $\B G_n$ from the spaces $\B G_{n-p}$ for $p > 0$.  The $n$th spectral sequence has $E^1$ page
\begin{align*}
& E^1_{p,q} \cong H_q(stab(\sigma_p);\bZ) \cong H_q(G_{n-p-1};\bZ), \\
  &E^1_{-1,q}  \cong H_q(G_{n};\bZ), 
 \end{align*}
  and $E^1_{p,q}=0$ for $p<-1$. 


The assumption that the complexes $W_n$ are highly connected implies that the spectral sequence converges to $0$ for $p+q \leq \frac{n-1}{2}$. The differential $$ d^1\colon  E^1_{0,i} = H_i(G_{n-1};\bZ) \longrightarrow E^1_{-1, i} = H_i(G_{n};\bZ)$$ is the map induced by the inclusion $G_{n-1} \hookrightarrow G_n$. Under the hypotheses of the theorem, we can argue by induction on $i$ that this map is an isomorphism (respectively, a surjection) in the desired range, to complete the proof of \autoref{QuillenStabillityTheorem}. 

In practice, given \autoref{QuillenStabillityTheorem}, the most difficult step in a proof of homological stability is usually the proof that the complexes $W_n$ are highly connected. 

In recent years, the argument that we just outlined has been axiomatized by Randal-Williams and Wahl  \cite{RWW}  and Krannich  \cite{Krannich}  to give  a very general framework to prove homological stability results, including  (co)homology  with twisted abelian and polynomial coefficients.  Another axiomatization is due to Hepworth  \cite{Hepworth}.

\subsection{An example: the braid group $\gB_n$} 

Let $\bD^2$ be the closed disk. Fix $n$ marked points in its interior and a distinguished point $*\in\partial \bD^2$.  Associated to the braid group $\gB_n$  is an $(n-1)$-dimensional simplicial complex $W_n$ called the \emph{arc complex}  which we define combinatorially.
\begin{itemize}
\item {\bf vertices:} $W_n$ has a vertex for each isotopy class of embedded arcs in $\bD^2$  joining $*$ with one of the marked points. 
 \item {\bf $p$-simplices:}  A set of $(p+1)$ vertices spans a $p$-simplex if the corresponding isotopy classes can be represented by arcs that are pairwise disjoint except at their starting point $*$. 
  \end{itemize}
  \begin{figure}[h]
\raggedleft

    \includegraphics[width=\columnwidth]{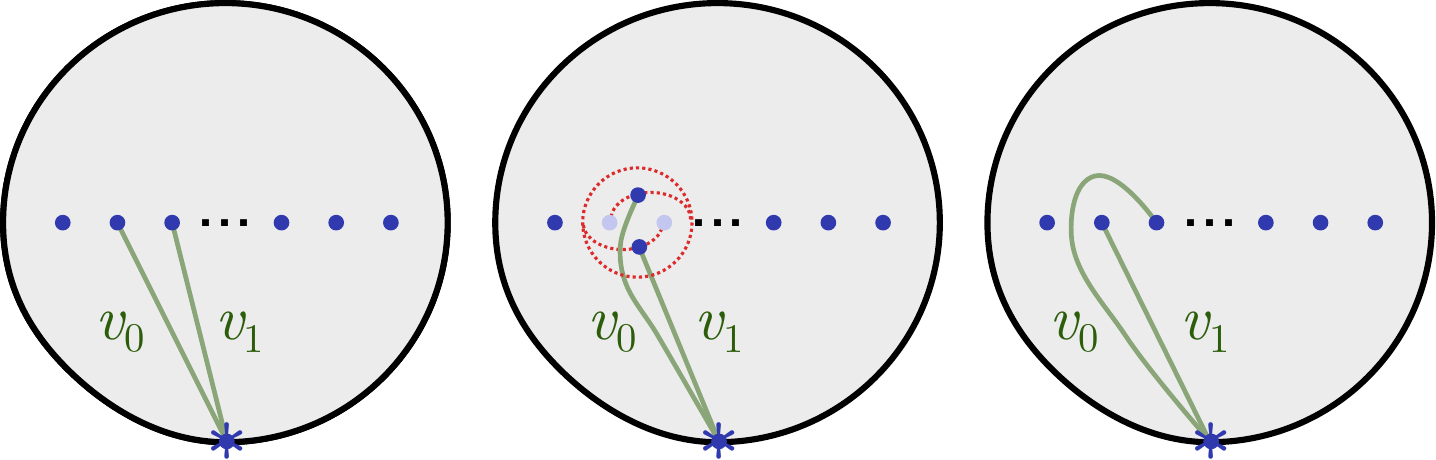}
\caption{  The action of  $\sigma_2\in\gB_n$ on a $1$-simplex $\{v_0, v_1\}$  of the arc complex  $W_n$.} \label{Arcs}

    \end{figure}  

Hatcher and Wahl \cite{HW}  proved that $W_n$ is $\left(\frac{n-2}{2}\right)$-connected (it is in fact contractible, see \cite{Damiolini}).

The braid group $\gB_n$ is isomorphic to the group $\Mod^n(\bD^2)$ of isotopy classes of diffeomorphisms of the closed disk that stabilize the set of marked points and restrict to the identity on  $\partial \bD^2$. Thus  $\gB_n$ has an action on $W_n$ that is simplicial and satisfies conditions $(i)$-$(v)$. See \autoref{Arcs}. \autoref{QuillenStabillityTheorem} gives a modern proof of homological stability for $\gB_n$ (Theorem \ref{HomStabilityBraids}), a result originally due to Arnold. 






\section{Representation stability}

 
\subsection{Configuration spaces  revisited}

\begin{figure}[h!]
\raggedleft
\centering
\includegraphics[scale=.6]{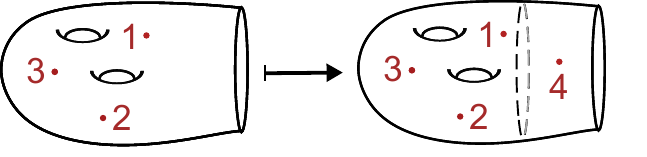} 
\caption{Stabilization map  $F_3(M)\rightarrow F_4(M)$} 
\label{OrderedStabilizationClean} 
\end{figure}

Let us address \autoref{QHomology} for the ordered configuration spaces $\{ F_n(M)\}_n$ when $M$ is the interior of a compact connected manifold with nonempty boundary. As with the unordered configuration spaces, given a choice of boundary component, we can define a stabilization map $F_n(M) \to F_{n+1}(M)$ that continuously introduces a new particle `at infinity'. See \autoref{OrderedStabilizationClean}.     

This suggests the question: for a fixed manifold $M$, do the spaces $\{ F_n(M)\}_n$ satisfy homological stability? The answer is, in contrast to $\{C_n(M)\}_n$, they do not, as we will verify directly. 
   
\begin{figure}[h!]
\raggedleft
\centering
\includegraphics[scale=.25]{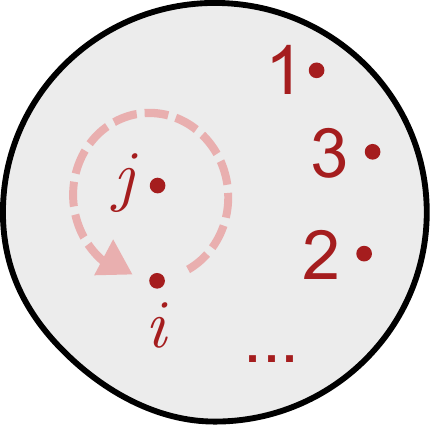} 
\caption{The homology class $\alpha_{ij} \in H_1(F_n(\C))$ } 
\label{H1FnD} 
\end{figure}

Let $M=\C$, so the homology $H_1(F_n(\C);\bZ)$ in degree 1 is the abelianization of the pure braid group $\gP_n$.  Artin's presentation implies that $\gP_n^{ab} \cong \Z^{{n \choose 2}}$ is free abelian on the images $\alpha_{ij}$ of the ${n \choose 2}$ generators $T_{ij}$ of  \autoref{PureBraidGenerator}. Viewed as a homology class in  $F_n(\C)$, we can represent $\alpha_{ij}$ by the loop illustrated in \autoref{H1FnD}. 
Hence, rank$(H_1(\gP_n;\bZ))$ grows quadratically in $n$, and homological stability fails. 


 Church and Farb \cite{CF},
 however, proposed a new paradigm for stability in spaces like the ordered configuration spaces $F_n(M)$ of a manifold $M$. Because (co)homology is functorial, the $S_n$-action  on $F_n(M)$ induces an action of $S_n$ on the (co)homology groups. Even though the (co)homology does not stabilize as a sequence of abelian groups, they proposed, it does stabilize as a sequence of $S_n$-representations. 
  
  There are several ways to formalize the idea of stability for a sequence of $S_n$--representations. One way, which was initially the primary focus of Church and Farb, is to consider the multiplicities of irreducible representations in the rational (co)homology groups. Suppose $V$ is a finite-dimensional rational $S_n$-representation. Because $S_n$ is a finite group, $V$ is \emph{semisimple}: it decomposes as a direct sum of {irreducible} subrepresentations. The multiplicities of the irreducible components  are uniquely defined and determine $V$ up to isomorphism.
  
  The irreducible rational $S_n$-representations are classified, and are in canonical bijection with partitions of $n$. A \emph{partition} $\lambda$ of a positive integer $n$ is a set of positive integers (called the \emph{parts} of $\lambda$) that sum to $n$. It is traditionally encoded by a \emph{Young diagram}, a collection of $n$ boxes arranged into rows of decreasing lengths equal to the parts of $\lambda$. For example, the Young diagram $\Y{3,2}$ corresponds to the partition $3+2$ of $5$. If $\lambda$ is a partition of $n$ (equivalently, a Young diagram of size $n$), we write $V_{\lambda}$ to denote the irreducible $S_n$-representation associated to $\lambda$. 
  
\begin{figure}[h!]
\centering
\includegraphics[scale=.8]{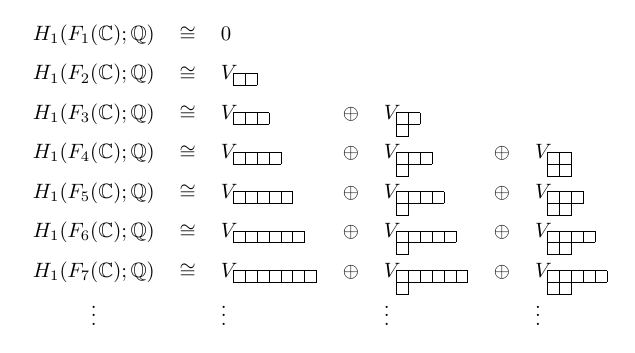}  
\caption{The decomposition of the homology groups $H_1(F_{n}(\C); \Q)$ for some small values of $n$.} 
\label{H1FCDecomp} 
\end{figure}

  Church and Farb observed a pattern in the homology of  $F_{n}(\C)$, which we illustrate in \autoref{H1FCDecomp} in homological degree $1$.


%

For $n \geq 4k$, we can recover the decomposition of $H_k(F_n(\C); \Q)$ into irreducible components simply by taking the decomposition of $H_k(F_{n-1}(\C); \Q)$ and adding a single box to the top row of each Young diagram. They showed that this pattern holds for all $k$, and Church \cite{Church} later proved that it holds for the cohomology groups $H^k(F_n(M); \Q)$ of the ordered configuration space of a connected oriented manifold of finite type. 

  Church, Farb, and others observed the same patterns in the (co)homology of a number of other families of groups and spaces. 
  These results raise the question, 
  \begin{question} \label{QStructure} 
  What underlying structure is responsible for these patterns?
  \end{question} Church, Ellenberg, Farb, Nagpal, Putman, and Sam answered this question by developing an algebraic framework that brought their work into a broader field, now called the field of {\it representation stability}. See, for example,  \cite{CF, CEF, CEFN, Putman,  CE, PutmanSam}. Other pioneers of the field, who approached it from different perspectives, include Sam, Snowden, Djament, Pirashvili, Vespa, Gan, and Li. Some selected references are \cite{Pirashvili, DjamentVespa-LinearGroups, SS-IntroTCAs, Snowden,  SS-Patterns, GanLi-EI, SS-GLEquivariant, Djament-Finiteness, SS-Grobner,  DjamentVespa-WeakPoly}. 
   
      





\subsection{$\FI$-modules} 

The key to answering Question \ref{QStructure} 
is the concept of an $\FI$-module. The theory of $\FI$--modules gives a conceptual framework that explains the ubiquity of the patterns observed in so many naturally-arising sequences of $S_n$-representations, and it also provides algebraic machinery to prove stronger results with streamlined arguments.

\begin{definition} 
Let $\FI$ be the category whose objects are finite sets (including $\varnothing$), and whose morphisms are all injective maps. Given a commutative ring $R$ (typically $\Z$ or $\Q$), an \emph{$\FI$-module} $V$ over $R$ is a functor from $\FI$ to the category of $R$-modules. 
\end{definition} 

\begin{figure}[h!]
\centering
\includegraphics[scale=3]{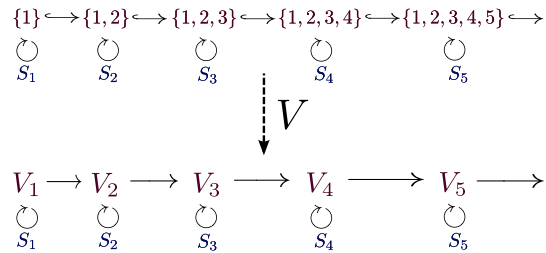}  
\caption{An $\FI$-module $V$ } 
\label{FIModule} 
\end{figure}

To describe an $\FI$-module $V$, it is enough to consider the ``standard" finite sets in $\FI$, $$[0]=\varnothing \quad \text{ and } \quad [n]= \{1, 2, \ldots, n\}.$$ For $n\geq 0$, we write $V_n$ to denote the image of $V$ on $[n]$. The endomorphisms of $[n]$ in $\FI$ are the symmetric group $S_n$, so $V_n$ is an $S_n$-representation. The data of an $\FI$-module $V$ is determined by the sequence of $S_n$-representations $\{V_n\}_n$, along with $S_n$-equivariant maps $\iota_n\colon  V_n \to V_{n+1}$ induced by the inclusion $[n] \hookrightarrow [n+1]$.  \autoref{FIModule} gives a schematic.

We refer to (the morphisms of) the category $\FI$ \emph{acting on} an $\FI$-module $V$ in the same sense that a ring $R$ acts on an $R$-module.

We encourage the reader to verify that the following sequences of $S_n$-representations form $\FI$-modules. 
\begin{itemize}
\item $V_n = \Q$ the trivial $S_n$--representations, \\
$\iota_n$ the identity maps.
\item $V_n = \Q^n$, $S_n$ permutes the standard basis, \\ $\iota_n\colon  \Q^n \cong (\Q^n \times \{0\}) \hookrightarrow \Q^{n+1}$.
\item $V_n = \Q[x_1, \ldots, x_n]$ the polynomial algebra with $S_n$ permuting the variables, $\iota_n$ the inclusion.
\end{itemize} 
Applying any endofunctor of $R$-modules to an $\FI$-module will produce another $\FI$-module, so we can construct more examples by (say)  taking tensor products or exterior powers of any of the above. 

We leave it as an exercise to the reader to verify that the following sequences of $S_n$-representations do {\bf not} form an $\FI$-module. A hint to this exercise: first verify that if  $\sigma \in S_n$ fixes the letters $\{1, 2, \ldots m\}$, then $\sigma$ must act trivially on the image of $V_m$ in $V_n$ under the map induced by the inclusion $[m] \subseteq [n]$. 
\begin{itemize}
\item $V_n = \Q$ the alternating representation,\\ i.e. $\sigma \cdot v = (-1)^{sgn(\sigma)} v \text{ for } v \in \Q$,\\  $\iota_n$ the identity map.
\item$V_n = \Q[S_n]$ the regular representation, \\ $\iota_n$ induced by the inclusion $S_n \subseteq S_{n+1}$. 
\end{itemize}

\begin{figure}[h!]
\raggedleft
\centering
 \includegraphics[width=\columnwidth]{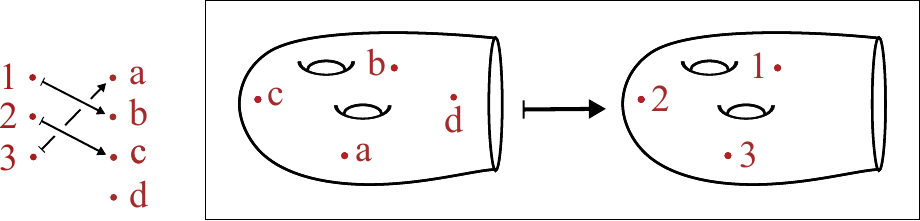}  
\caption{An $\FI$ morphism and its contravariant action  on the configuration spaces  $\{F_n(M)\}_n$} 
\label{CoFIModStabilization} 
\end{figure}

Importantly for present purposes, the (co)homology groups of ordered configuration spaces form $\FI$-modules in many cases. If $M$ is any space, there is a contravariant action of $\FI$ on its ordered configuration spaces by continuous maps. If we view a point in $F_n(M)$ as an embedding $\rho\colon  [n] \to M$, then an $\FI$ morphism $f\colon  [m]\to [n]$ acts by precomposition,
\begin{align*} 
f^*\colon  F_n(M) & \longrightarrow F_m(M) \\ 
\rho & \longmapsto \rho \circ f.
\end{align*}  See \autoref{CoFIModStabilization}. 

Composing this $\FI$ action with the (contravariant) cohomology functor gives a {\bf covariant} action of $\FI$ on the cohomology groups $\{H^k(F_n(M))\}_n$. 

\begin{figure}[h!]
\raggedleft
\centering
 \includegraphics[width=\columnwidth]{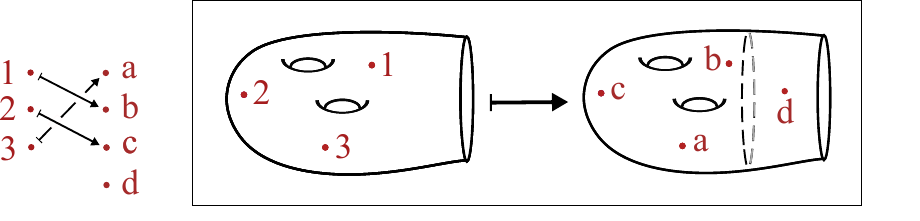} 
\caption{An $\FI$ morphism and its covariant action  on the configuration spaces  $\{F_n(M)\}_n$} 
\label{FIModStabilization} 
\end{figure}

To obtain a covariant action of $\FI$ on $\{F_n(M)\}_n$, we need additional assumptions on the space $M$. Let $M$ be the interior of a compact manifold of dimension at least 2 with nonempty boundary. Consider an $\FI$ morphism $f\colon  [m]\to [n]$ and a configuration in $F_m(M)$. We relabel particles by their image under $f$, and apply the stabilization map of Section \ref{HomStabConf} to introduce any particles not in $f([m])$ in a neighbourhood of a distinguished boundary component. See \autoref{FIModStabilization}. 

This  action of $\FI$ is only functorial up to homotopy, but this suffices to induce a well-defined $\FI$-module structure on the sequence of homology groups $\{H_k(F_n(M))\}_n$.

Modules over the category $\FI$ behave in many ways like modules over a ring (technically, they are an abelian category). We define a map of $\FI$-modules $V \to W$ to be a natural transformation, that is, a sequence of maps $V_n \to W_n$ that commute with the $\FI$ morphisms. The kernels and images of these maps themselves form $\FI$-modules, and we can define operations like tensor products and direct sums in a natural way. This structure allows us to import many of the standard tools from commutative and homological algebra to the study of $\FI$-modules. 

\begin{figure}[h!]
\raggedleft
\centering
\includegraphics[scale=3.1]{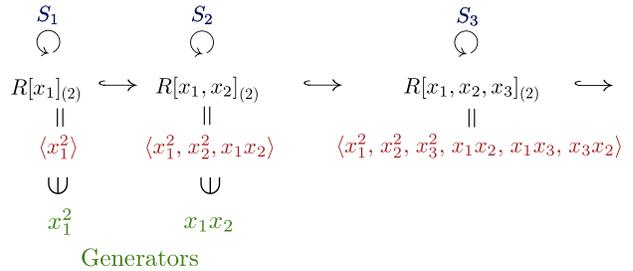} 
\caption{A finite generating set for the $\FI$-module $R[x_1, \ldots x_n]_{(2)}$ } 
\label{PolyDeg2Generators} 
\end{figure}

 Church, Ellenberg, and Farb \cite{CEF} showed the answer to  Question \ref{QStructure} is that the sequences in question are $\FI$-modules that are {finitely generated}. 
\begin{definition}  Let $V$ be an $\FI$-module. A subset $S \subseteq \bigsqcup_{n \geq 0} V_n $  \emph{generates} $V$ if the images of $S$ under the $\FI$ morphisms span $V_n$ for all $n \geq 0$. Equivalently, the smallest $\FI$-submodule of $V$ containing $S$ is $V$ itself. The $\FI$-module $V$ is \emph{finitely generated} in \emph{degree $\leq d$} if there is a finite subset of elements $S \subseteq \bigsqcup_{ n \leq d} V_n$ that generates $V$. 
\end{definition} 
For example, consider the $\FI$-module $V$ over a ring $R$ such that $V_n = R[x_1, \ldots, x_n]_{(d)}$ is the submodule of homogeneous degree-$d$ polynomials in $n$ variables, $S_n$ acts by permuting the variables, and $\iota_n\colon V_n \to V_{n+1}$ is the inclusion map. We encourage the reader to verify that $V$ is finitely generated in degree $\leq d$. \autoref{PolyDeg2Generators} shows a finite generating set when  $d=2$. 

\begin{figure}[h!]
\raggedleft
\centering
\includegraphics[scale=.27]{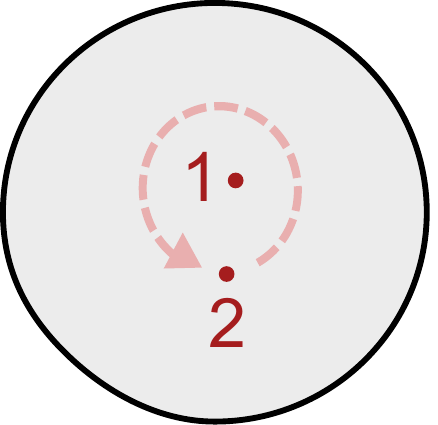} 
\caption{The homology class $\alpha_{1,2} \in H_1(F_2(\C))$ generates the $\FI$-module $\{H_1(F_n(\C);\Q)\}_n$ } 
\label{H1F2D} 
\end{figure}

Another example: from our description of the groups $\{H_1(F_n(\C); \Q)\}_n$ in \autoref{H1FnD}, we see that this $\FI$-module is generated by the single element $\alpha_{1,2} \in H_1(F_2(\C);\Q)$ shown in \autoref{H1F2D}. 
Arnold's description of the homology groups of $F_n(\C)$ \cite{ArnoldRing} makes it straightforward to verify finite generation of $\{H_k(F_n(\C);\Q)\}_n$ in every degree $k$.

Church, Ellenberg  and Farb \cite{CEF}, and independently Snowden \cite{Snowden} proved that $\FI$-modules over $\bQ$ satisfy a \emph{Noetherian} property: submodules of finitely generated modules are themselves always finitely generated.  Using this result, Church--Ellenberg--Farb proved that, if $V$ is a finitely generated $\FI$-module, then the sequence $\{V_n\}_n$ of $S_n$--representations stabilizes in several senses. 

\begin{theorem}[Church--Ellenberg--Farb \cite{CEF}] \label{RepStability} Let $V$ be an $\FI$-module over $\Q$,  finitely generated in degree $\leq d$. The following hold. 
\begin{itemize} 
\item {\bf Finite generation.} For $n \geq d$, $$S_{n+1} \cdot \iota_n(V_n)\text{ \  spans \ }V_{n+1}.$$ 
\item {\bf Polynomial growth.} There is a polynomial in $n$ of degree $\leq d$ that agrees with the dimension $\dim_{\Q}(V_n)$ for all $n$ sufficiently large. 
\item {\bf Multiplicity stability.} For all $n \geq 2d$ the decomposition of $V_n$ into irreducible constituents stabilizes (in the sense illustrated in \autoref{H1FCDecomp}). 
\item {\bf Character polynomials.} The character of $V_n$ is independent of $n$ for all $n \geq 2d$.
\end{itemize} 
\end{theorem}
The characters of $V$ are in fact eventually equal to a \emph{character polynomial}, independent of $n$; see \cite[Section 3.3]{CEF}. 

The answer of  \autoref{QHomology}  for the family $\{F_n(M)\}_n$ is then given by the following result.






\begin{theorem}[Church \cite{Church}; Church--Ellenberg--Farb \cite{CEF};  Miller--Wilson \cite{MW}]  \label{RepStabilityConfig}  Let $M$ be  the interior of a compact connected smooth manifold of dimension at least 2 with nonempty boundary. In each degree $k$ the homology and cohomology of ordered configuration spaces $\{F_n(M)\}_n$ of $M$ are finitely generated $\FI$-modules. In particular, the rational (co)homology groups stabilize in the sense of \autoref{RepStability}.
\end{theorem} 

\begin{figure}[h!]
\raggedleft
\centering
\includegraphics[scale=.7]{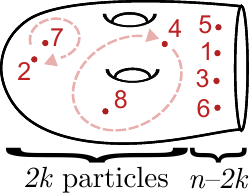} 
\caption{A homology class in the image of $H_k(F_{2k}(M);\bZ)$.} 
\label{ImageOfOrderedStabilization} 
\end{figure}

Heuristically, \autoref{RepStabilityConfig} states that the homology of $F_n(M)$ is spanned by classes of the form shown in \autoref{ImageOfOrderedStabilization}. 
%


From the $S_n$-covering relationship (\autoref{CnMQuotient})  it follows that  $\dim H^k(C_n(M);\mathbb{Q})$  is equal to the multiplicity of the trivial representation  in $H^k(F_n(M);\mathbb{Q})$. Hence Theorem \ref{RepStabilityConfig}  implies classical cohomological stability with $\mathbb{Q}$-coefficients for unordered configuration spaces $\{ C_n(M)\}_n$. Church \cite{Church} used representation stability techniques to prove rational (co)homological stability results for the unordered configuration spaces $\{ C_n(M)\}_n$ even in the case that $M$ is a closed manifold, so the isomorphisms are not necessarily induced by natural stabilization maps. See also Randal-Williams \cite{RWConf}.


\subsection{Other instances of representation stability} 

The definition of a finitely generated $\FI$-module makes sense for representations over the integers or other coefficients, even in situations where the representations are not semisimple and multiplicity stability is not well-defined.  Moreover, this approach readily generalizes to analogous categories that encode actions by families of groups other than the symmetric groups. Some examples that have been studied are the classical Weyl groups, certain wreath products, various linear groups, and products or decorated variants of $\FI$. The term ``representation stability'' now refers to algebraic finiteness results (like finite generation or presentation degree)  for a module over one of these categories. 
For further reading on representation stability, see the introductory notes and articles \cite{FarbICM, Wilson-FINotes, Snowden-FINotes, Sam-Notices}.

The (co)homology of several families of groups and moduli spaces exhibit representation stability.\\[-.5em]

\noindent{\bf Generalized ordered configuration spaces and pure braid groups.}
There is a large and growing body of work on representation stability for the homology of configuration spaces: improving stable ranges, studying configuration spaces of broader classes of topological spaces, or studying alternate stabilization maps. See for example \cite{E-WG, casto,    Petersen,HR, lutgehetmann,  KupersHomotopy, CMNR, bahran, MW, Ramos, Alpert}.

Other families generalizing the pure braid groups also have representation stable cohomology groups, including the  {pure virtual braid groups}, the {pure  flat braid groups}, the {pure cactus groups}, and the {group of pure string motions}  \cite{WilsonS, LeeVirtual, JRMD}. \\[-.5em]

\noindent{\bf Pure mapping class groups and moduli spaces of surfaces with marked points}.
Given a set of $n$ labelled marked points in a surface $\Sigma$, the {\it mapping class group} $\Mod^n(\Sigma)$ is   the group of isotopy classes of (orientation-preserving  if $\Sigma$ is orientable) diffeomorphisms of $\Sigma$ that fix $\partial\Sigma$ and stabilize the set of marked points. The {\it pure mapping class group} $\PMod^n(\Sigma)$ is the subgroup that fixes the marked points pointwise. These groups  also generalize the braid groups  since   $\Mod^n(\bD^2)\cong\gB_n$ and $\PMod^n(\bD^2)\cong\gP_n$. There is a short exact sequence
$$ 1 \to \PMod^n(\Sigma) \to \Mod^n(\Sigma) \to S_n \to 1$$
that defines an action of $S_n$ on the (co)homology of $\PMod^n(\Sigma)$. Hatcher and Wahl \cite{HW} proved that the sequence $\{\Mod^n(\Sigma)\}_n$ satisfies homological stability. Jim\'enez Rolland  \cite{JR, JR1, JR2} proved that the  groups $H^k(\PMod^n(\Sigma);\bZ)$ assemble to a finitely generated $\FI$-module.

 For $g\geq 2$ the {\it moduli space $\mathcal{M}_{g,n}$ of Riemann surfaces of genus $g$ with $n$ marked points} is a rational model of the classifying space $\B\PMod^n(\Sigma_g)$, and the symmetric group $S_n$ acts on $\mathcal{M}_{g,n}$ by permuting the $n$ marked points.  Hence,   the sequence $\{H^k(\mathcal{M}_{g,n};\bQ)\}_n$ of $S_n$-representations  stabilizes in the sense of  Theorem \ref{RepStability}.  
 
 In contrast,  for fixed genus $g$ the cohomology groups $H^k
(\overline{\mathcal{M}}_{g,n}; \bQ)$ of the {\it Deligne-Mumford compactification of $\mathcal{M}_{g,n}$} can grow exponentially in $n$. Thus these sequence cannot be finitely generated as  $\FI$-modules.  Tosteson \cite{Tosteson} proved, however, that the sequences $\{H^k (\overline{\mathcal{M}}_{g,n}; \bQ)\}_n$ are subquotients of finitely generated $\FS^{op}$-modules, where  $\FS^{op}$ is the opposite category of the category of finite sets and surjective maps.  From this he deduced constraints on the growth rate and on the irreducible $S_n$-representations that occur.\\[-.5em]

\noindent{\bf Flag varieties.} Let $\gG^\mathcal{W}_n$  be a semisimple complex Lie group of type $A_{n-1}$, $B_n$, $C_n$, or $D_n$, with Weyl group $\mathcal{W}_n$ and
$\gB^{\mathcal{W}}_n$ a Borel subgroup. The space $\gG^\mathcal{W}_n/\gB^\mathcal{W}_n$  is called a {\it generalized flag variety}. 
Representation stability of these cohomology groups (as $S_n$- or $\mathcal{W}_n$-representations) has been studied by Church--Ellenberg--Farb \cite{CEF}, Wilson \cite{Wilson1}, and others.\\[-.5em]

\noindent{\bf Complements of arrangements.} The cohomology of hyperplane complements associated to certain reflection groups $\mathcal{W}_n$  (and  their toric and elliptic analogues) stabilizes as a sequence of $\mathcal{W}_n$-representations by the work of Wilson \cite{Wilson2} and Bibby \cite{Bibby}. Representation stability holds for the cohomology of more general linear subspace arrangements  with a wider class of groups actions; see Gadish \cite{Gadish}.\\[-.5em] 

\noindent{\bf Congruence subgroups.}  
Let $K$ be a commutative ring and $I \subseteq K$ a proper two-sided ideal. The  {\it level $I$ congruence subgroups} $\GL_n(K, I)$ of $\GL_n(K)$ are defined to be the kernel of
the ``reduction modulo $I$'' map $\GL_n(K)\rightarrow \GL_n(K/I)$. 
Representation stability of the sequence of homology groups $\{H_k(\GL_n(K,I);\bZ)\}_n$ (as $S_n$ or $\GL_n(K/I)$-representations) has been studied by  Gan-- Li \cite{GanLi},  Putman \cite{Putman},  Putman--Sam \cite{PutmanSam}, Church--Ellenberg--Farb--Napgal \cite{CEFN}, Miller--Patzt--Wilson \cite{MPW-CentralStability}, Miller--Nagpal--Patzt \cite{MNP} and others.










\section{Current research directions}

Work continues on proving (co)homological stability for new families or new coefficients systems, improving stable ranges, and computing the stable and unstable (co)homology for families known to stabilize. 

Recently Galatius, Kupers  and  Randal-Williams  \cite{GKRW-Ek, GKRW-GL, GKRW-MCG} identified and  proved a new kind of  stabilization result, which  they describe by the slogan ``the failure of homological stability is itself stable''.  They defined homological-degree-shifting stabilization maps and use them to prove {\it secondary homological stability} for the homology of mapping class groups and general linear groups.  Himes \cite{Himes} studied secondary stability for unordered configuration spaces. Miller--Patzt--Petersen \cite{MPP-polynomial} studied stability with polynomial coefficient systems. Miller--Wilson  \cite{MW},  Bibby--Gadish \cite{BibbyGadish}, Ho \cite{QHo}, and Wawrykow \cite{Wawrykow} studied representation-theoretic analogues of secondary stability for ordered configuration spaces. 

For a more in-depth introduction to homological stability and these current research directions, we recommend Kupers' minicourse notes \cite{Kupers} and references therein. 





\section{Acknowledgments} We thank Omar Antol\'in Camarena, Jeremy Miller, and Nicholas Wawrykow for useful feedback on a draft of this article. We thank our referees for their extensive comments.  The authors are grateful to Benson Farb and Tom Church for introducing them to the field of representation stability, and to many of the ideas in this survey. Rita Jim\'enez Rolland is grateful for the financial support  by the CONACYT grant Ciencia Frontera
CF-2019/217392. Jennifer Wilson is grateful for the support of NSF grant DMS-1906123. 

\bibliographystyle{alpha}
\bibliography{StabilityRefs.bib} 


\end{document}